# RENEWALS FOR EXPONENTIALLY INCREASING LIFETIMES, WITH AN APPLICATION TO DIGITAL SEARCH TREES

By Florian Dennert and Rudolf Grübel

*Universität Hannover*

We show that the number of renewals up to time $t$ exhibits distributional fluctuations as $t \to \infty$ if the underlying lifetimes increase at an exponential rate in a distributional sense. This provides a probabilistic explanation for the asymptotics of insertion depth in random trees generated by a bit-comparison strategy from uniform input; we also obtain a representation for the resulting family of limit laws along subsequences. Our approach can also be used to obtain rates of convergence.

**1. Introduction.** Let $(Y_k)_{k \in \mathbb{N}}$ be a sequence of independent, nonnegative random variables and let $(S_n)_{n \in \mathbb{N}_0}$,

$$ S_0 := 0, \qquad S_n := \sum_{k=1}^{n} Y_k \qquad \text{for all } n \in \mathbb{N}, $$

be the associated sequence of partial sums. Regarding the $Y_k$'s as successive lifetimes and $S_n$ as the time of the $n$th renewal, we interpret

$$ N_t := \sup\{n \in \mathbb{N}_0 : S_n \le t\} $$

as the number of renewals up to and including time $t$; $(N_t)_{t \ge 0}$ is the renewal process. Standard renewal theory assumes that the $Y_k$'s all have the same distribution, in which case $N_t$, appropriately rescaled, is asymptotically normal as $t \to \infty$. For this result, and for renewal theory in general, we refer the reader to Section XI in [3].

In this note we consider exponentially increasing lifetimes. We show that in such a case the distribution of $N_t$ does not converge and that asymptotic distributional fluctuations appear (Section 2). Such fluctuations occur frequently in the analysis of algorithms. The renewal theoretic framework provides a probabilistic view of this phenomenon in connection with digital search trees (Section 3). We also indicate how our approach can be used to obtain rates of convergence (Section 4).

**2. Renewals for increasing lifetimes.** We assume that the lifetimes increase exponentially with rate $\alpha$, where $\alpha > 1$ is fixed throughout the sequel, in the sense that

$$(1) \qquad \alpha^{-k} Y_k \to_{\text{distr}} Y_\infty \quad \text{and} \quad \alpha^{-k} E Y_k \to E Y_\infty$$

for some random variable $Y_\infty$ and as $k \to \infty$. Here "$\to_{\text{distr}}$" denotes convergence in distribution, so that the first part of (1) means that

$$\lim_{k \to \infty} E f(\alpha^{-k} Y_k) = E f(Y_\infty)$$

for all bounded continuous functions $f : \mathbb{R} \to \mathbb{R}$. Below we will use the fact that in order to prove $X_n \to_{\text{distr}} X$ it is sufficient to show that $E f(X_n) \to E f(X)$ holds for all bounded and *uniformly* continuous functions. For details and a general treatment of convergence in distribution we refer the reader to [1]. Of course, only the distribution $\mu = \mathcal{L}(Y_\infty)$ of $Y_\infty$ matters, so we will occasionally write $\alpha^{-k} Y_k \to_{\text{distr}} \mu$ instead. Finally, throughout this note a condition involving moments is meant to imply that these moments are finite.

An important role will be played by

$$S_\infty := \sum_{k=0}^{\infty} \alpha^{-k} Y_{\infty, k},$$

where $(Y_{\infty, k})_{k \in \mathbb{N}_0}$ is a sequence of independent and identically distributed random variables with $\mathcal{L}(Y_{\infty, 0}) = \mathcal{L}(Y_\infty)$, $Y_\infty$ as in (1). From $E Y_\infty < \infty$ we obtain $E S_\infty = \alpha (\alpha - 1)^{-1} E Y_\infty < \infty$ and therefore $P(S_\infty < \infty) = 1$; moreover, we then also have that $\sum_{k=0}^{n} \alpha^{-k} Y_{\infty, k}$ converges almost surely and hence in distribution to $S_\infty$ as $n \to \infty$. We will also assume that $\mathcal{L}(Y_\infty)$ has no atoms, that is,

$$(2) \qquad P(Y_\infty = y) = 0 \qquad \text{for all } y \in \mathbb{R}_+.$$

Finally, it is an elementary analytic fact that, for a sequence $(x_n)_{n \in \mathbb{N}}$ of real numbers with limit $x \in \mathbb{R}$,

$$(3) \qquad \lim_{n \to \infty} \sum_{k=0}^{n-1} \alpha^{-k} x_{n-k} = \frac{\alpha x}{\alpha - 1}.$$

The following lemma can be regarded as a random version of (3).

LEMMA 1. *If* (1) *and* (2) *are satisfied, then* $\alpha^{-n} S_n \to_{\text{distr}} S_\infty$ *as* $n \to \infty$, *and* $P(S_\infty = y) = 0$ *for all* $y \in \mathbb{R}$.

Proof. Suppose that $(U_k)_{k \in \mathbb{N}}$ is a sequence of independent random variables on some probability space $(\Omega, \mathcal{A}, P)$, all uniformly distributed on the unit interval. Let $F_k$ be the distribution function of $Y_k$, $F$ the distribution function of $Y_\infty$. We use a variant of the quantile construction:

$$\tilde{Y}_k := F_k^{-1}(U_k), \qquad \tilde{Y}_{\infty,k} := F^{-1}(U_k) \qquad \text{for all } k \in \mathbb{N}.$$

We then have $\mathcal{L}(\tilde{Y}_1, \dots, \tilde{Y}_n) = \mathcal{L}(Y_1, \dots, Y_n)$ for all $n \in \mathbb{N}$, which implies

$$\mathcal{L}(\alpha^{-n} S_n) = \mathcal{L}(\alpha^{-n} \tilde{S}_n) \qquad \text{with } \tilde{S}_n := \sum_{k=1}^{n} \tilde{Y}_k.$$

Using $\alpha^{-n} \tilde{S}_n = \sum_{k=0}^{n-1} \alpha^{-k} (\alpha^{-(n-k)} \tilde{Y}_{n-k})$ we obtain

$$(4) \qquad E \left| \alpha^{-n} \tilde{S}_n - \sum_{k=0}^{n-1} \alpha^{-k} \tilde{Y}_{\infty,n-k} \right| \leq \sum_{k=0}^{n-1} \alpha^{-k} E |\alpha^{-(n-k)} \tilde{Y}_{n-k} - \tilde{Y}_{\infty,n-k}|.$$

With $Y_k' := F_k^{-1}(U_1)$ and $Y_\infty' := F^{-1}(U_1)$ we have

$$(5) \qquad E |\alpha^{-k} \tilde{Y}_k - \tilde{Y}_{\infty,k}| = E |\alpha^{-k} Y_k' - Y_\infty'|.$$

From (1) it follows that $\alpha^{-k} Y_k' \to_{\mathrm{distr}} Y_\infty'$ and $E \alpha^{-k} Y_k' \to E Y_\infty'$. Because of $Y_k', Y_\infty' \geq 0$ Theorem 5.4 in [1] applies and gives the $L^1$-convergence of $\alpha^{-k} Y_k'$ to $Y_\infty'$, that is, $E |\alpha^{-k} Y_k' - Y_\infty'| \to 0$ as $k \to \infty$. Using this together with (3), (4) and (5) we obtain

$$(6) \qquad \lim_{n \to \infty} E \left| \alpha^{-n} \tilde{S}_n - \sum_{k=0}^{n-1} \alpha^{-k} \tilde{Y}_{\infty,n-k} \right| = 0.$$

Now let $f : \mathbb{R} \to \mathbb{R}$ be bounded and uniformly continuous. We have

$$|E f(\alpha^{-n} S_n) - E f(S_\infty)| = \left| E f(\alpha^{-n} \tilde{S}_n) - E f\left( \sum_{k=0}^{n-1} \alpha^{-k} \tilde{Y}_{\infty,n-k} \right) \right.$$

$$\left. + E f\left( \sum_{k=0}^{n-1} \alpha^{-k} \tilde{Y}_{\infty,k} \right) - E f\left( \sum_{k=0}^{\infty} \alpha^{-k} \tilde{Y}_{\infty,k} \right) \right|$$

$$\leq E \left| f(\alpha^{-n} \tilde{S}_n) - f\left( \sum_{k=0}^{n-1} \alpha^{-k} \tilde{Y}_{\infty,n-k} \right) \right|$$

$$+ E \left| f\left( \sum_{k=0}^{\infty} \alpha^{-k} \tilde{Y}_{\infty,k} \right) - f\left( \sum_{k=0}^{n-1} \alpha^{-k} \tilde{Y}_{\infty,k} \right) \right|.$$

For the first integral on the right-hand side we use (6), for the second an elementary estimate shows that the difference between the arguments of $f$ converges to 0 in probability. In both cases we now use uniform continuity when the arguments of $f$ are close to each other and boundedness otherwise. This leads to

$$\lim_{n \to \infty} Ef(\alpha^{-n} S_n) = Ef(S_\infty),$$

which gives the convergence in distribution. The statement about the atoms of $S_\infty$ follows from (2) and the fact that $S_\infty$ is equal in distribution to $Y_\infty + \alpha^{-1} S_\infty$ with $Y_\infty$ and $S_\infty$ independent. □

The above proof is based on classical weak convergence arguments. An alternative proof can be obtained via the Wasserstein distance

$$d_W(\mu, \nu) = \inf\{E|X - Y| : \mathcal{L}(X) = \mu, \mathcal{L}(Y) = \nu\},$$

its known relation to weak convergence and convergence of the first moments, and the same variant of the quantile construction, which in this context is known as the comonotone coupling.

We write $\lfloor x \rfloor$ for the greatest integer less than or equal to $x$ and $\{x\}$ for the fractional part of $x \in \mathbb{R}$.

THEOREM 2.  *Suppose that* (1) *and* (2) *are satisfied and let*

$$(7) \qquad Q_\eta := \mathcal{L}(\lfloor -\log_\alpha S_\infty + \eta \rfloor), \qquad 0 \le \eta \le 1.$$

*If* $(t_n)_{n \in \mathbb{N}}$ *is a sequence of real numbers with* $t_n \to \infty$ *and such that* $\{\log_\alpha t_n\} \to \eta$ *for some* $\eta \in [0, 1]$, *then*

$$N_{t_n} - \lfloor \log_\alpha t_n \rfloor \to_{\text{distr}} Q_\eta \qquad \text{as } n \to \infty.$$

PROOF.  We use the abbreviations $k_n := \lfloor \log_\alpha t_n \rfloor$ and $\eta_n := \{\log_\alpha t_n\}$. In particular, $\log_\alpha t_n = k_n + \eta_n$. Further, let $Z_\infty := -\log_\alpha S_\infty$. By a standard renewal theoretic argument,

$$P(N_t = j) = P(S_j \le t) - P(S_{j+1} \le t) \qquad \text{for all } t \ge 0, j \in \mathbb{N}_0,$$

hence

$$\begin{aligned}
P(N_{t_n} - k_n = j) &= P(S_{k_n + j} \le t_n) - P(S_{k_n + j + 1} \le t_n) \\
&= P(-\log_\alpha(\alpha^{-k_n - j} S_{k_n + j}) + \eta_n \ge j) \\
&\quad - P(-\log_\alpha(\alpha^{-k_n - j - 1} S_{k_n + j + 1}) + \eta_n \ge j + 1) \\
&\to P(\lfloor Z_\infty + \eta \rfloor = j) \qquad \text{as } n \to \infty,
\end{aligned}$$

where in the last step Lemma 1 and three general facts about convergence in distribution were used: First, the continuous mapping theorem, which implies that $-\log_\alpha(\alpha^{-m} S_m) \to_{\text{distr}} -\log_\alpha S_\infty$ as $m \to \infty$; secondly, the interplay with convergence in probability, see Theorem 4.1 in [1], which yields

$-\log_\alpha(\alpha^{-n}S_n) + \eta_n \to_{\text{distr}} -\log_\alpha S_\infty + \eta$ as $n \to \infty$; finally, that $\mathcal{L}(S_\infty)$ and therefore also $\mathcal{L}(-\log_\alpha S_\infty + \eta)$ assign probability 0 to single points and that this implies

$$\lim_{n\to\infty} P(-\log_\alpha(\alpha^{-n}S_n) + \eta_n \geq z) = P(-\log_\alpha S_\infty + \eta \geq z) \qquad \text{for all } z \in \mathbb{R}.$$

$\square$

A structural consequence of the representation (7) is the $\to_{\text{distr}}$-continuity of $\eta \mapsto Q_\eta$ on the open unit interval; at $\eta = 0$ this function is right continuous, at $\eta = 1$ it is left continuous. The extreme members are translates of each other in the sense that $Q_0(\{j\}) = Q_1(\{j+1\})$ for all $j \in \mathbb{Z}$.

The total variation distance $d_{\text{TV}}$ of probability measures is defined by

$$d_{\text{TV}}(\mu, \nu) := \sup_B |\mu(B) - \nu(B)|,$$

for $\mu, \nu$ concentrated on $\mathbb{Z}$ this can be written as

$$(8) \qquad d_{\text{TV}}(\mu, \nu) = \tfrac{1}{2} \sum_{j\in\mathbb{Z}} |\mu(\{j\}) - \nu(\{j\})|.$$

For a sequence of probability measures that are concentrated on a fixed countable set Scheffé's lemma implies that weak convergence is equivalent to convergence in total variation distance, hence (7) can be rewritten as

$$\lim_{n\to\infty} d_{\text{TV}}(\mathcal{L}(N_{t_n} - \lfloor \log_\alpha t_n \rfloor), Q_{\{\log_\alpha t_n\}}) = 0.$$

Because of the continuity of $[0,1] \ni \eta \mapsto Q_\eta$ this in turn leads to a statement that avoids the use of subsequences,

$$(9) \qquad \lim_{t\to\infty} d_{\text{TV}}(\mathcal{L}(N_t - \lfloor \log_\alpha t \rfloor), Q_{\{\log_\alpha t\}}) = 0.$$

In Section 4 we will investigate the rate of convergence in (9) in a particular case.

## 3. An application to digital search trees.
The nodes of a (rooted, directed) binary tree can be represented by finite strings of 0's and 1's if we interpret 0 as a move to the left and 1 as a move to the right. The length of the string is the depth (or level) of the node it represents, the root node corresponds to the empty string and has level 0. The sequence $(T_n)_{n\in\mathbb{N}}$ associated with a sequence $(x_n)_{n\in\mathbb{N}}$ of numbers from the unit interval by the DST (digital search tree) algorithm is obtained as follows: For $T_1$, we put $x_1$ into the root node. If $x_1, \ldots, x_n$ have been stored into $T_n$ then the position of $x_{n+1}$ is determined by traveling through the tree with the direction given by the binary expansion of $x_{n+1}$ until an empty node has been found. This algorithm and its properties are discussed in the standard texts of the area, for example, [8, 10, 11]. As an example we consider the first ten numbers given in [8],

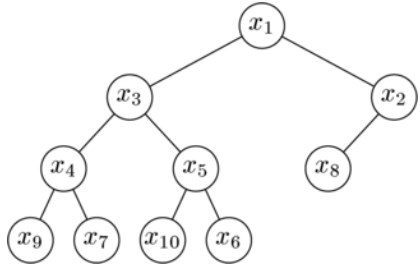

Appendix A, $(\sqrt{2}, \sqrt{3}, \sqrt{5}, \sqrt{10}, \sqrt[3]{2}, \sqrt[3]{3}, \sqrt[4]{2}, \log 2, \log 3, \log 10)$. Let $x_i$ be the fractional part of the $i$th entry, $1 \le i \le 10$; the relevant first four bits of the respective binary expansions are given by $(0110, 1011, 0011, 0010, 0100, 0111, 0011, 1011, 0001, 0100)$. This leads to the binary tree given in Figure 1.

Consider now the sequence $(T_n)_{n \in \mathbb{N}_0}$ of random trees that the DST algorithm associates with a sequence $(U_n)_{n \in \mathbb{N}}$ of independent random variables, where we assume that the $U_n$'s are uniformly distributed on the unit interval and that $T_0$ is the empty tree. Let $X_n(\theta)$ be the depth of the first free node of $T_n$ along the path determined by a sequence $\theta \in \{0, 1\}^{\mathbb{N}}$. Such a $\theta$ defines a family of nested intervals of length $2^{-k}$, $k = 1, 2, 3, \ldots$, and it is easy to see that $(X_n(\theta))_{n \in \mathbb{N}_0}$ is a Markov chain with $X_0(\theta) \equiv 0$ and transition probabilities $p_{k,k+1} = 1 - p_{k,k} = 2^{-k}$ for all $k \in \mathbb{N}_0$. Conditioning on the value of $U_{n+1}$ we see that the distribution of $X_n(\theta)$ is the same as the distribution of $Z_{n+1}$, the insertion depth of $U_{n+1}$. This quantity is known as "unsuccessful search" in the literature on the analysis of algorithms. [Of course, this distributional equality does not hold for the joint distributions: $n \mapsto X_n(\theta)$ is increasing, $n \mapsto Z_{n+1}$ is not.] For example, the next number in Knuth's list is $x_{11} = 1/\log 2$, the binary expansion of the fractional part $\{x_{11}\}$ begins with 011100 and hence $x_{11}$ would be inserted at level 4 as the right child of $x_6$.

The Markov chain $(X_n(\theta))_{n \in \mathbb{N}_0}$ is of the simple birth type and can therefore be described by its respective holding times $Y_1, Y_2, Y_3, \ldots$ in the states $k = 0, 1, 2, \ldots$. These are independent, and $Y_k$ has a geometric distribution with parameter $p_{k-1,k}$, that is, for all $k \in \mathbb{N}$,

$$P(Y_k = j) = (1 - 2^{-k+1})^{j-1} 2^{-k+1} \qquad \text{for all } j \in \mathbb{N}.$$

Here we interpret the case $k = 1$ as $Y_1$, the holding time in 0, being constant and equal to 1. As a result of its simple stochastic dynamics, $(X_n(\theta))_{n \in \mathbb{N}_0}$ is equal to the renewal process $N$ associated with the sequence $(Y_k)_{k \in \mathbb{N}}$, observed at discrete times, that is, $(X_n(\theta))_{n \in \mathbb{N}_0} = (N_n)_{n \in \mathbb{N}_0}$. It is easy to see that for this sequence $(Y_k)_{k \in \mathbb{N}}$ of lifetimes conditions (1) and (2) are satisfied and that $\mathcal{L}(Y_\infty) = \text{Exp}(2)$, with $\text{Exp}(\lambda)$ the exponential distribution with parameter $\lambda$ (and mean $1/\lambda$). Hence Theorem 2 can be applied: If the sequence $(n(m))_{m \in \mathbb{N}} \subset \mathbb{N}$ is such that $n(m) \to \infty$ and $\{\log_2 n(m)\} \to \eta$ as $m \to \infty$, then

$$(10) \qquad X_{n(m)}(\theta) - \lfloor \log_2 n(m) \rfloor \to_{\text{distr}} Q_\eta.$$

Here $Q_\eta$, $0 \le \eta \le 1$, is the distribution of $\lfloor -\log_2 S_\infty + \eta \rfloor$, $S_\infty := \sum_{k=0}^{\infty} 2^{-k} Y_{\infty,k}$ and $Y_{\infty,k}$, $k \in \mathbb{N}_0$, are independent and identically distributed with $\mathcal{L}(Y_{\infty,1}) = \text{Exp}(2)$. Alternatively, we can write $S_\infty := \sum_{k=1}^{\infty} \tilde{Y}_k$ with $\tilde{Y}_k$, $k \in \mathbb{N}$, again independent and $\mathcal{L}(\tilde{Y}_k) = \text{Exp}(2^k)$ for all $k \in \mathbb{N}$.

The explicit representation of the family of limit distributions on the basis of the convolution product of the distributions $\text{Exp}(2^k)$, $k \in \mathbb{N}$, can be used to obtain a series expansion for the distribution functions associated with $Q_\eta$, $0 \le \eta \le 1$. For this, we start with a partial fraction expansion: For all $n \in \mathbb{N}$ and all $z \in \mathbb{C}$ with $|\Re(z)| < 2$,

$$(11) \qquad \prod_{k=1}^{n} (1 - 2^{-k} z)^{-1} = \sum_{k=1}^{n} a_{n,k} (1 - 2^{-k} z)^{-1},$$

where $a_{n,k} := \prod_{j=1}^{k-1} (1 - 2^j)^{-1} \prod_{j=1}^{n-k} (1 - 2^{-j})^{-1}$. Reading (11) as an equality relating characteristic functions we obtain

$$(12) \qquad \text{Exp}(2^1) \star \text{Exp}(2^2) \star \cdots \star \text{Exp}(2^n) = \sum_{k=1}^{n} a_{n,k} \text{Exp}(2^k).$$

Note, however, that the right-hand side in (12) is not the usual mixture of probability distributions as the coefficients alternate in sign. With

$$a_k := b \prod_{j=1}^{k-1} (1 - 2^j)^{-1}, \qquad b := \prod_{j=1}^{\infty} (1 - 2^{-j})^{-1},$$

letting $n \to \infty$ in (12) leads to $\mathcal{L}(S_\infty) = \sum_{k=1}^{\infty} a_k \text{Exp}(2^k)$, so that

$$
\begin{aligned}
(13) \qquad Q_\eta((-\infty, x]) &= P(\lfloor -\log_2(S_\infty) + \eta \rfloor \le x) \\
&= P(S_\infty > 2^{\eta - 1 - x}) \\
&= \sum_{k=1}^{\infty} a_k \exp(-2^{k + \eta - 1 - x}) \qquad \text{for all } x \in \mathbb{Z}.
\end{aligned}
$$

This representation of the limiting distribution functions as an alternating series has already been obtained by Louchard [9] in the context of digital search trees and by Flajolet [4] in the context of approximate counting; see also Section 6.4 in [10] and Section 6.3 in [8] for related results. These authors use a completely different approach, more analytic in flavor and relying on combinatorial identities due to Euler.

Our main point here, however, is not a rederivation of (13) but the representation of the family $\{Q_\eta : 0 \leq \eta < 1\}$ in terms of one particular random variable, which is first shifted by $\eta$ and then discretized. This representation can, for example, be used to obtain information about the tail behavior of the limit distributions. Janson [7] notes that (13) by itself would only give an exponential rate of decrease for the tail probabilities, he then provides an analytic argument that improves this to a superexponential rate by showing that the associated Fourier transform is an entire function. Using the representation $S_\infty = \sum_{k=1}^{\infty} 2^{-k} Z_k$ with $Z_k$ independent and $\mathcal{L}(Z_k) = \mathrm{Exp}(1)$ together with the fact that $\mathrm{Exp}(1)$ has a density bounded by 1, we get

$$P(S_\infty \leq 2^{-j}) \leq P(Z_1 \leq 2^{-j+1}) P(Z_2 \leq 2^{-j+2}) \cdots P(Z_{j-1} \leq 2^{-1})$$
$$\leq 2^{-j+1} 2^{-j+2} \cdots 2^{-1}$$
$$= 2^{-j(j-1)/2}$$

for all $j \in \mathbb{N}$. Because of $Q_\eta([k, \infty)) \leq P(S_\infty \leq 2^{-k+1})$ for all $k \in \mathbb{N}$, $k \geq 2$, this leads to

$$Q_\eta([x, \infty)) = o(\exp(-\rho x^2)) \qquad \text{as } x \to \infty, \text{ for all } \rho < (\log 2)/2.$$

The fact that a representation by discretization is possible in many situations where fluctuations were first found by calculation seems to belong to the folklore of the subject, at least in simple instances such as the asymptotic distributional behavior of the maximum of a sample from a geometric distribution. The geometric case together with some renewal theoretic techniques (for identically distributed lifetimes) was used in [5] to obtain results of the above type for von Neumann addition. In [2] a discretization representation occurs on the level of stochastic processes, leading to a probabilistic approach to fluctuation phenomena in the context of multiplicities of the maximum in a random sample from a discrete distribution. In a recent paper, Janson [7] studies the effects of discretizing random variables and the resulting distributional fluctuations and gives a range of interesting examples. Of course, the explanation for periodicities can be, and indeed often is, quite different and mechanisms other than discretization may be responsible; see, for example, [6] and the references given there.

**4. Rates of convergence.** The renewal theoretic approach can also be used to obtain rates of convergence. We sketch one of the possibilities, for a particular choice of distances, and give details for the DST situation from Section 3. Let, for $t > 0$, $k(t) := \lfloor \log_\alpha t \rfloor$ and $\eta(t) := \{\log_\alpha t\}$.

The Kolmogorov–Smirnov distance of two probability measures $\mu$ and $\nu$ on the real line is defined by

$$d_{\mathrm{KS}}(\mu, \nu) := \sup_{x \in \mathbb{R}} |\mu((-\infty, x]) - \nu((-\infty, x])|.$$

If $X$ and $Y$ are real random variables, then we abbreviate $d_{KS}(\mathcal{L}(X), \mathcal{L}(Y))$ to $d_{KS}(X, Y)$; if $F$ and $G$ are the associated distribution functions, then $d_{KS}(X, Y) = \|F - G\|_\infty$, where the supremum norm for general bounded functions $f : \mathbb{R} \to \mathbb{R}$ is given by $\|f\|_\infty := \sup_{x \in \mathbb{R}} |f(x)|$. The Kolmogorov–Smirnov distance is obviously invariant under strictly monotone transformations. For example,

$$d_{KS}(\alpha X + \beta, \alpha Y + \beta) = d_{KS}(X, Y) \qquad \text{for all } \alpha, \beta \in \mathbb{R}, \alpha \neq 0,$$

and for $X, Y > 0$,

$$d_{KS}(X, Y) = d_{KS}(\log X, \log Y).$$

With the notation as in the proof of Theorem 2,

$$\begin{aligned}
|P(N_t - k(t) = j) &- P(\lfloor -\log_\alpha(S_\infty) + \eta(t) \rfloor = j)| \\
&\leq |P(-\log_\alpha(\alpha^{-k(t)-j} S_{k(t)+j}) + \eta(t) \geq j) - P(-\log_\alpha(S_\infty) + \eta(t) \geq j)| \\
&\quad + |P(-\log_\alpha(\alpha^{-k(t)-j-1} S_{k(t)+j+1}) + \eta(t) \geq j + 1) \\
&\qquad\qquad - P(-\log_\alpha(S_\infty) + \eta(t) \geq j + 1)|.
\end{aligned}$$

With the auxiliary quantities

$$Z_t := \lfloor -\log_\alpha(S_\infty) + \eta(t) \rfloor, \qquad \phi(m) := d_{KS}(\alpha^{-m} S_m, S_\infty)$$

and the above properties of the Kolmogorov–Smirnov distance this leads to

$$(14) \qquad |P(N_t - k(t) = j) - P(Z_t = j)| \leq \phi(k(t) + j) + \phi(k(t) + j + 1).$$

It is often possible to obtain an upper bound for negative $j$, say $j \leq -k(t)/2$, directly. In such cases the above elementary renewal theoretic argument leads to a bound for the $\|\cdot\|_\infty$-distance between the probability mass functions of $N_t - k(t)$ and $Z_t$, for example; note that the latter variable has distribution $Q_{\eta(t)}$ where $Q_\eta$, $0 \leq \eta \leq 1$, is the set of limit distributions along subsequences that appears in Theorem 2.

The above argument covers the step from $(\alpha^{-m} S_m)_{m \in \mathbb{N}}$ to $(N_t)_{t \geq 0}$. However, in an application the starting point will usually be the convergence of the scaled lifetimes in (1), which means that we also need an analogue for Lemma 1 that gives rates of convergence.

We carry this out in the specific context of digital search trees. The following general bounds will turn out to be useful: If $X$ has density $f_X$ and if $P(|Y| \leq c) = 1$, then

$$(15) \qquad\qquad d_{KS}(X, X + Y) \leq c \|f_X\|_\infty.$$

Indeed: For all $z \in \mathbb{R}$, $P(X \leq z - c) \leq P(X + Y \leq z) \leq P(X \leq z + c)$, so that

$$\begin{aligned}
|P(X + Y \leq z) &- P(X \leq z)| \\
&\leq \max\{P(X \leq z + c) - P(X \leq z), P(X \leq z) - P(X \leq z - c)\},
\end{aligned}$$

and, of course, $P(X \in (a, b]) \leq (b - a)\|f_X\|_\infty$. This bound can easily be generalized to

(16)     $$d_{\mathrm{KS}}(X, X + Y) \leq c\|f_X\|_\infty + P(|Y| > c) \qquad \text{for all } c > 0,$$

where we still assume that $X$ has density $f_X$, but $Y$ may be arbitrary. Note that $X$ and $Y$ need not be independent in (15) and (16). If they are independent then it is easy to show that

(17)     $$d_{\mathrm{KS}}(X, X + Y) \leq \|f_X\|_\infty E|Y|.$$

In (17) boundedness of $Y$ is not needed but the bound obviously makes sense only if $Y$ has finite first moment. Finally, in connection with density bounds the interplay with convolution is of interest: We have $\|f \star g\|_\infty \leq \|f\|_\infty$ for all probability densities $f, g$. For example, if a sum of independent random variables contains a summand with distribution $\mathrm{Exp}(\lambda)$, then the density of the sum is bounded by $\lambda$.

LEMMA 3.   *With $(Y_k)_{k\in\mathbb{N}}$ and $S_\infty$ as in Section* 3,
$$d_{\mathrm{KS}}(2^{-n}S_n, S_\infty) = O(n2^{-n}).$$

PROOF.   Let $(Z_k)_{k\in\mathbb{N}}$ be a sequence of independent random variables, all exponentially distributed with parameter 1. Then $S_\infty$ is equal in distribution to $\sum_{k=1}^\infty 2^{-k}Z_k$. We recall that the $k$th lifetime $Y_k$ has a geometric distribution with parameter $2^{-k+1}$. On the basis of $(Z_k)_{k\in\mathbb{N}}$ we define a sequence $(\tilde{Y}_k)_{k\in\mathbb{N}}$ by $\tilde{Y}_k := \lfloor \alpha_k Z_k \rfloor + 1$ for all $k \in \mathbb{N}$, with

$$\alpha_1 := 0, \qquad \alpha_k := (-\log(1 - 2^{-k+1}))^{-1} \qquad \text{for } k > 1.$$

It is easy to check that

$$(\tilde{Y}_k)_{k\in\mathbb{N}} =_{\mathrm{distr}} (Y_k)_{k\in\mathbb{N}}, \qquad 2^{-n}\sum_{k=1}^n 2^{k-1}Z_k =_{\mathrm{distr}} \sum_{k=1}^n 2^{-k}Z_k.$$

Hence, with $\phi(n)$ denoting the $d_{\mathrm{KS}}$-distance of $2^{-n}S_n$ and $S_\infty$,

$$\phi(n) \leq \phi_1(n) + \phi_2(n) + \phi_3(n) \qquad \text{for all } n \in \mathbb{N},$$

with $\phi_1, \phi_2, \phi_3$ defined by

$$\phi_1(n) := d_{\mathrm{KS}}\left(2^{-n}\sum_{k=1}^n \tilde{Y}_k, 2^{-n}\sum_{k=1}^n \alpha_k Z_k\right),$$

$$\phi_2(n) := d_{\mathrm{KS}}\left(2^{-n}\sum_{k=1}^n \alpha_k Z_k, 2^{-n}\sum_{k=1}^n 2^{k-1}Z_k\right),$$

$$\phi_3(n) := d_{\mathrm{KS}}\left(\sum_{k=1}^n 2^{-k}Z_k, \sum_{k=1}^\infty 2^{-k}Z_k\right).$$

For the random variables in $\phi_1$ we have

$$V_n \leq 2^{-n} \sum_{k=1}^n \tilde{Y}_k \leq V_n + n2^{-n} \qquad \text{with } V_n := 2^{-n} \sum_{k=1}^n \alpha_k Z_k.$$

It is easy to show that the densities of $V_n$, $n \in \mathbb{N}$, can be uniformly bounded for all $n$ by some finite constant $C_1$, hence (15) implies that $\phi_1(n) \leq C_1 n2^{-n}$ for all $n \in \mathbb{N}$.

The elementary bounds

$$\frac{1}{x} - 1 \leq -\frac{1}{\log(1-x)} \leq \frac{1}{x} \qquad \text{for } 0 < x \leq \frac{1}{2}$$

together with $\alpha_1 = 0$ imply $\sup_{k \in \mathbb{N}} |\alpha_k - 2^{k-1}| = 1$, hence we have

$$\left| 2^{-n} \sum_{k=1}^n \alpha_k Z_k - 2^{-n} \sum_{k=1}^n 2^{k-1} Z_k \right| \leq 2^{-n} \sum_{k=1}^n Z_k.$$

The familiar combination of Markov's inequality and moment generating functions gives

$$P\left( \sum_{k=1}^n Z_k \geq (1 + \kappa)n \right) = O(2^{-n})$$

if $\kappa$ is chosen large enough, so that we can use (16) with $c = c(n) = (1 + \kappa)n2^{-n}$ to obtain that $\phi_2(n) \leq C_2 n2^{-n}$ for all $n \in \mathbb{N}$, for some finite constant $C_2$.

For $\phi_3$ finally we use (17): For the densities of the finite sums we again have a finite uniform bound for all $n$, and

$$E\left| \sum_{k=n+1}^\infty 2^{-k} Z_k \right| = \sum_{k=n+1}^\infty 2^{-k} E Z_k = 2^{-n},$$

so that $\phi_3(n) \leq C_3 2^{-n}$ for all $n \in \mathbb{N}$ with some $C_3 < \infty$. Putting these together we arrive at

$$\phi(n) \leq Cn2^{-n} \qquad \text{for all } n \in \mathbb{N}$$

with some finite constant $C$. $\quad \square$

In the application under consideration we obtain a rate of convergence result with respect to the total variation distance, which is stronger than a result for the supremum norm distance of the corresponding probability mass functions that we mentioned in connection with (14).

THEOREM 4. *With $(X_n(\theta))_{n \in \mathbb{N}}$ and $Q_\eta$ as in Section 3,*

$$d_{\mathrm{TV}}(\mathcal{L}(X_n(\theta) - \lfloor \log_2 n \rfloor), Q_{\{\log_2 n\}}) = o(n^{-\gamma}) \qquad \text{for all } \gamma < 1.$$

PROOF.   We use the abbreviations $k(n) := \lfloor \log_2 n \rfloor$ and $\eta(n) := \{\log_2 n\}$. Let $\gamma < 1$ be given and choose $\varepsilon > 0$ such that $\varepsilon < 1 - \gamma$. Lemma 3 together with (14) gives

$$\sum_{j \geq -\varepsilon k(n)} |P(N_n - k(n) = j) - Q_{\eta(n)}(\{j\})| \leq C \sum_{j \geq (1-\varepsilon)k(n)} j2^{-j}$$

for all $n \in \mathbb{N}$ with some finite constant $C$. Our choice of $\varepsilon$ implies that the upper bound has the desired rate $o(n^{-\gamma})$.

For the remaining part of the infinite sum in (8) we replace the absolute difference of the probabilities by their sum, which means that it is now enough to show that

$$(18) \qquad\qquad P(N_n \leq (1 - \varepsilon)k(n)) = o(n^{-\gamma}),$$

$$(19) \qquad\qquad P(-\log_2(S_\infty) \leq -\varepsilon k(n) + 1) = o(n^{-\gamma}).$$

Here we have used that $Q_\eta$ is the distribution of $\lfloor -\log_2(S_\infty) + \eta \rfloor$. It is easy to show that the moment generating function for $S_\infty$ exists in a neighborhood of 0, hence

$$(20) \qquad\qquad P(S_\infty > x) = o(e^{-\kappa x}) \qquad \text{for all } x > 0$$

with some $\kappa > 0$. Straightforward manipulations show that (20) implies (19); indeed, the probability converges faster to 0 than any negative power of $n$. Using once again the relation between the number of renewals and the partial sums of the lifetimes we further obtain, with $m(n, \varepsilon) := \lfloor (1 - \varepsilon)k(n) \rfloor$,

$$P(N_n \leq (1 - \varepsilon)k(n)) \leq P(S_{m(n,\varepsilon)} \geq n)$$
$$= P(2^{-m(n,\varepsilon)} S_{m(n,\varepsilon)} \geq n2^{-m(n,\varepsilon)})$$
$$\leq d_{\mathrm{KS}}(2^{-m(n,\varepsilon)} S_{m(n,\varepsilon)}, S_\infty) + P(S_\infty \geq n2^{-m(n,\varepsilon)}).$$

For the Kolmogorov–Smirnov distance we use Lemma 3, for the tail of $S_\infty$ the desired rate follows with (20). This gives (18) and hence completes the proof.  □

## REFERENCES

[1] BILLINGSLEY, P. (1968). *Convergence of Probability Measures.* Wiley, New York. MR0233396

[2] BRUSS, F. TH. and GRÜBEL, R. (2003). On the multiplicity of the maximum in a discrete random sample. *Ann. Appl. Probab.* **13** 1252–1263. MR2023876

[3] FELLER, W. (1971). *An Introduction to Probability Theory and Its Applications* **II**, 2nd ed. Wiley, New York. MR0270403

[4] FLAJOLET, PH. (1985). Approximate counting: A detailed analysis. *BIT* **25** 113–134. MR0785808

[5] GRÜBEL, R. and REIMERS, A. (2001). On the number of iterations required by von Neumann addition. *Theor. Inform. Appl.* **35** 187–206. MR1862462

[6] Janson, S. (2004). Functional limit theorems for multitype branching processes and generalized Pólya urns. *Stochastic Process. Appl.* **110** 177–245. MR2040966

[7] Janson, S. (2006). Rounding of continuous random variables and oscillatory asymptotics. *Ann. Probab.* **34** 1807–1826.

[8] Knuth, D. E. (1973). *The Art of Computer Programming* **3**. *Sorting and Searching*. Addison–Wesley, Reading, MA. MR0445948

[9] Louchard, G. (1987). Exact and asymptotic distributions in digital binary search trees. *Theor. Inform. Appl.* **21** 479–496. MR0928772

[10] Mahmoud, H. M. (1992). *Evolution of Random Search Trees.* Wiley, New York. MR1140708

[11] Sedgewick, R. and Flajolet, Ph. (1996). *An Introduction to the Analysis of Algorithms.* Addison–Wesley, Reading, MA.

Institut für Mathematische Stochastik
Universität Hannover
Postfach 60 09
D-30060 Hannover
Germany
E-mail: dennert@stochastik.uni-hannover.de
rgrubel@stochastik.uni-hannover.de